\documentclass[11pt]{article}
\newcommand{\newsection}[1]{\setcounter{equation}{0}\section{#1}}

\def\keywords{ \if@twocolumn
\section*{Keywords}
\else \small
\begin{center}
{ \bf Keywords\vspace{-.5em}\vspace{0pt}}
\end{center}
\center \fi}
\def\endkeywords{ \if@twocolumn\else\endcenter\fi}
\usepackage{epsfig}
\usepackage{graphicx}
\begin{document}
\title{\bf Hierarchical cascade model leading to $7$-th order initial
value problem}
\author{\thanks{
Department of Mathematics, University of the Punjab, Lahore 54590,
 Pakistan.   ~~~~~~~~~~~~~~~~~~~~~~~~~~~~~~~~~~~~Email: toghazala2003@yahoo.com} Ghazala Akram \thanks{
School of Mathematical Sciences, Queen Mary, University of London,
Mile End Road, London E1 4NS, UK. ~~Email:
 g.akram@qmul.ac.uk} ,
 ~Christian Beck \thanks{School of Mathematical Sciences, Queen Mary, University of London,
Mile End Road, London E1 4NS, UK.~~~Email:
 c.beck@qmul.ac.uk}}
\date{}
\maketitle
\begin{center}

\abstract{In turbulent flows, local velocity differences often obey
a cascade-like hierarchical dynamics, in the sense that local
velocity differences at a given scale $k$ are driven by
deterministic and random forces from the next-higher scale $k-1$.
Here we consider such a hierarchically coupled model with periodic
boundary conditions, and show that it leads to an $N$-th order
initial value problem, where $N$ is the number of cascade steps. We
deal in detail with the case $N=7$ and introduce a non-polynomial
spline method that solves the problem for arbitrary driving forces.
Several examples of driving forces are considered, and estimates of
the numerical precision of our method are given. We show how to
optimize the numerical method to obtain a truncation error of order
$O(h^5)$ rather than $O(h^2)$, where $h$ is the discretization
step.}

\end{center}

\section{Introduction}
Hierarchical dynamics arise quite commonly for complex systems that
consist of many subdynamics that drive each other in a selfsimilar
way. Often, the  dynamics at a given scale $k$ couples to the
dynamics at a higher scale $k-1$ in a simple way, and these
types of problems can sometimes be dealt with in an analytic way. A
typical example are cascade models in turbulence (see, e.g.
\cite{cascade,prl,sawford} and references therein), which are a useful tool to
characterize the selfsimilar features of turbulent flows at high
Reynolds numbers.

In this paper we consider a very simple hierarchical model that can
be physically interpreted as representing velocity differences at
different scales that are driven by deterministic and stochastic
forces from the next higher scale, in a medium with viscosity and a
given driving force at the top scale. We show that the problem
reduces to the solution of an initial value problem associated with
a differential equation of order $N$, where $N$ is the number of
cascade steps. As an example, we deal in detail with the case $N=7$,
which corresponds to typical values of cascade sizes observed in
turbulent flows. There are three standard approaches to solve
initial value problems numerically, the finite difference method,
the finite element method and the spline approximation methods. We
introduce a non-polynomial spline method for the numerical solution
of this initial value problem, and estimate the precision of our
numerical treatment. Higher-order boundary value problems are
effectively solved using non-polynomial spline methods 
\cite{Sg074,
Sg076,Sg071,Sg07,Sg08}. It turns
out there is an optimum choice of the interpolation parameters where
our methods yields best possible results (5th order rather than 2nd
order in the discretization step). Some analytically solvable
examples of driving forces are dealt with as examples, for
illustration of our general method.

\section{Cascade model}

Let us consider a model of damped particles in a viscous medium
that are driven by rapidly fluctuating forces.
Suppose there are $N$ such particles. We denote the velocity
of each particle as $y^{(k)}(t)$, $k=1,2,\ldots , N$.
A very simple, uncoupled model would be that the particles
are damped by a linear friction force, which is proportional to velocity,
and a rapidly fluctuating driving force
$L^{(k)}(t)$ that is independent of velocity:
\begin{equation}
\dot{y}^{(k)}=-\Gamma y^{(k)}+ L^{(k)}(t)
\end{equation}
$\Gamma >0$ denotes the friction constant.
If $L^{(k)}(t)$ is Gaussian white noise, then this model
just leads to the Ornstein-Uhlenbeck process, performed independently
by each particle \cite{vKa}.

Here, however, we want to modify this model to a more interesting
interacting dynamics.
First of all, we allow $L^{(k)}(t)$ to be any time-dependent driving force,
and assume that it is differentiable. Next, we construct a coupled
hierarchical model, by replacing $y^{(k)}$ on the right-hand side
of the above equation by the
nearest neighbour $y^{(k+1)}$. We may physically
interprete $y^{(k)}$ as a local velocity difference in a turbulent
flow at spatial scale $r=2^{-k}$.
The physical interpretation is that the change of velocity due to
friction forces at a given scale $k$ is proportional to the velocity at the
next-smaller scale $k+1$. This reminds us of the fact that in cascade
models of turbulence energy dissipates from larger scales
down to smaller scales. In our hierarchical model the actual dissipation
force at spatial level $k$ is proportional to the velocity difference
at the next smaller
scale $k+1$:
\begin{equation}
\dot{y}^{(k)}= -\Gamma y^{(k+1)}+L^{(k)} (t) \;\;\;\;k=1,\ldots ,N \label{here}
\end{equation}
We now show that this model, with a cascade of size $N$, leads to
an initial value problem of the form
\begin{equation}
\frac{ \partial^n}{\partial t^n} y= (-1)^N \Gamma^N y +g(t), \label{london}
\end{equation}
where $y(t)=y^{(1)}(t)$ is the velocity difference at the
top of the cascade ($k=1$), $g(t)$ is a driving force at the top of the
cascade, and we have assumed periodic boundary conditions at the top
and bottom of the cascade, i.e. $y^{(n+1)}(t)=y^{(1)}(t)$.

To derive eq.~(\ref{london}), we differentiate eq.~(\ref{here}) to obtain
\begin{eqnarray}
\ddot{y}^{(k)}&=&-\Gamma \dot{y}^{(k+1)}+\dot{L}^{(k)} (t) \\
&=& +\Gamma^2 y^{(k+2)}-\Gamma L^{(k+1)}(t)+ \dot{L}^{(k)}(t),
\end{eqnarray}
where in the last step we used eq.~(\ref{here}) with $k$ replaced by  $k+1$.
Further differentiation yields
\begin{eqnarray}
\frac{\partial^3}{\partial t^3} {y}^{(k)} &=& \Gamma^2 \dot{y}^{(k+2)}
+\ddot{L}^{(k)}(t)-\Gamma \dot{L}^{(k+1)}(t) \\
&=&-\Gamma^3 y^{(k+3)}+\Gamma^2 L^{(k+2)}(t)-\Gamma \dot{L}^{(k+1)}(t)+\ddot{L}^{(k)} (t).
\end{eqnarray}
Finally, for a cascade with $N$ steps one arrives at
\begin{equation}
\frac{\partial^N}{\partial t^N} y^{(k)}(t)=(-1)^n \Gamma^N y^{(k+N)}+g(t),
\end{equation}
where $g(t)$ is a linear combination of derivatives of the driving
forces at the various scales:
\begin{equation}
g(t)=(-\Gamma)^{N-1} L^{(k+N-1)} (t)+(-\Gamma)^{N-2} \dot{L}^{(k+N-2)}
+ \ldots
\end{equation}
Our implemented periodic boundary
condition of the cascade
$y^{(1)}(t)=y^{(N)}(t)$
simply means that at smallest scales the dynamics
should just be the same as at the largest scales,
which is a self-similarity assumption.
For $N$ odd and defining $y^{(1)}(t)=y(t)$ we arrive at
\begin{equation}
\frac{\partial^N}{\partial t^N}y+\Gamma^N y=g(t),
\end{equation}
which has to be supplemented by a set of $N$ initial
conditions, corresponding to initial velocities at the various scales.
This initial value problem will be solved in the next section.

\newsection{The initial value problem}
In the following, we choose as an example $N=7$. In turbulence simulations,
the driving forces at various scales of the turbulent flow are only known in numerical form,
and given the chaotic nature of the
forces it is important to implement high precision
numerical methods that optimize the numerical solution of driven velocity
fields within a given
cascading subdynamics. In the following, we introduce a
high-precision non-polynomial spline method for the
solution of the initial value problem
(IVP) that corresponds to our cascade model. \\\\
Consider the following seventh order initial value problem
\begin{equation}\label{e:1.1}
\left. \begin{array}{lcl}
    y^{(7)}(t)+f(t)y(t)&=&g(t),  \ \ \ \ t \ \in \  [a, \ b], \\
    y(a)  \ \ \ = u_{0},\ & &  y^{(1)}(a) \ \ \  =  u_{1}, \\
    y^{(2)}(a) = u_{2},\ & & y^{(3)}(a)  =   u_{3}, \\
    y^{(4)}(a) = u_{4}, \ & & y^{(5)}(a) = u_{5}, \\
    y^{(6)}(a) = u_{6}, \ & &
    \end{array} \right \}
\end{equation}
where $ u_{i} ; i= 0,1,\ldots,6$ are finite real constants while the
functions $f(t)$ and $g(t)$ are continuous on $[a, b]$. This is just
the dynamics derived in the previous section, provided the effective
friction constant $\Gamma$ depends on time $t$. The notation is
slightly different since now $y^{(i)}$ denotes the $i$-th derivative
with respect to $t$. For our physical application in terms of a
cascade-like model, we need to keep the functions $f(t)$ and $g(t)$
quite general since they are unknown in a turbulent flow. Our aim is
to provide a proper numerical method to provide a most accurate
solution of this IVP. It turns out that a non-polynomial spline
method is very useful in this context. For particular choices of the
interpolation parameters our method provides optimum results (error
term of order five rather than 2), as shown in the following
sections.

  \newsection{Nonpolynomial Spline Method}
  \indent To develop the spline approximation to the problem $(\ref{e:1.1})$,
  the interval $[a, \ b]$ is divided
  into $n$ equal subintervals, using the grid points $t_{i}  =  a +  ih $ ;
   $i=0,1,\ldots ,
  n$, where $h \ = \ (b-a)/n.$ \\
    Consider the following restriction $S_{i}$ of the approximate solution $S$ to each subinterval $[t_{i}, \
  t_{i+1}], \ i=0,1,\ldots , n-1$,
  \begin{eqnarray}
  S_{i}(t)& = &a_{i}\cos \omega(t-t_{i})+b_{i}\sin \omega(t-t_{i})+c_{i}(t-t_{i})^{6}
  +d_{i}(t-t_{i})^{5}+e_{i}(t-t_{i})^{4}
                 \nonumber \\
               & &  +
                p_{i}(t-t_{i})^{3}+q_{i}(t-t_{i})^{2}+r_{i}(t-t_{i})+v_{i}.
  \end{eqnarray}
     Let
  \begin{equation}
 \left. \begin{array}{rcl}
  y_{i} \ = \ S_{i}(t_{i}) & & \ \ \ \ m_{i} \ = \ S_{i}^{(1)}(t_{i}),
    \\
  M_{i} \ = \ S_{i}^{(2)}(t_{i}),  & & \ \ \ \ N_{i} \ = \
  S_{i}^{(4)}(t_{i}), \\
  U_{i} \ = \ S_{i}^{(7)}(t_{i}),  & &
  \end{array} \right \} \ i=0,1,\ldots,n.
  \end{equation}
  Following ~\cite{Sg076} and postulating continuous derivatives at knots, consistency relations
between the values of splines and their seventh order derivatives at
knots are obtained as
\begin{eqnarray} \label{e:2.8}
& \left(\alpha h^{7}U_{i-7}+\beta h^{7}U_{i-6}+\gamma
h^{7}U_{i-5}+\delta h^{7}U_{i-4}+\delta h^{7}U_{i-3}+\gamma
h^{7}U_{i-2}+\beta
h^{7}U_{i-1}+\alpha h^{7}U_{i} \right) \nonumber \\
= & \left[-120y_{i-7} +840 y_{i-6} -  2520 y_{i-5} +  4200 y_{i-4}-
4200 y_{i-3}
  +  2520 y_{i-2}  - 840  y_{i-1} +120 y_{i}  \right]; \nonumber \\
  & i=7,8,\ldots,n,
\end{eqnarray}
where
\begin{eqnarray}
\alpha & = & \left(\frac{120(\cos \theta-1)}{\theta^{7}\sin \theta}+
\frac{60}{ \theta^{5} \sin \theta}-\frac{5}{
\theta^{3}\sin\theta}+\frac{1}{6 \theta \sin\theta}\right),
\nonumber
\\
\beta & = & \left(\frac{600(1-\cos \theta)}{\theta^{7}\sin \theta}-
\frac{60(2 \cos \theta -3)}{ \theta^{5} \sin \theta}+\frac{5(2 \cos
\theta -9)}{ \theta^{3}\sin\theta}-\frac{(2
\cos \theta -57)}{6 \theta \sin\theta}\right), \nonumber \\
\gamma & = & \left(\frac{1080(\cos \theta-1)}{\theta^{7}\sin
\theta}+ \frac{180(2 \cos \theta +1)}{ \theta^{5} \sin
\theta}+\frac{45(2 \cos \theta +1)}{ \theta^{3}\sin\theta}-\frac{(38
\cos \theta -101)}{2 \theta \sin\theta}\right) \nonumber
\end{eqnarray}
and
\begin{eqnarray}
 \delta & = &
\left(\frac{600(1-\cos \theta)}{\theta^{7}\sin \theta}- \frac{60(4
\cos \theta +1)}{ \theta^{5} \sin \theta}-\frac{5(20 \cos \theta
-1)}{ \theta^{3}\sin\theta}-\frac{(604 \cos \theta -359)}{6 \theta
\sin\theta}\right). \nonumber
\end{eqnarray}
Here $\theta=\omega h$ is an arbitrary parameter. The relation
$(\ref{e:2.8})$ forms a system of $(n-6)$ linear equations in the
$(n)$ unknowns $(y_{i}, \ i=1,2,...,n) $, while $U_{i}$ is taken
from IVP
$(\ref{e:1.1})$ to be equal to $-f_{i}y_{i} + g_{i}$, $i=0,1,\ldots,n.$ \\
 Six further equations (end conditions) are required to obtain
a complete solution for the $y_{i}s$ appearing in eq.
$(\ref{e:2.8})$. These equations are calculated using the method of
undetermined coefficients \cite{Sg074} as follows
\begin{eqnarray}\label{e:3.1}
  U_{0} -10 U_{1} +  U_{4} & = & \frac{1}{h^{7}} \ \left[  \ \frac{512540}{27} y_{0}
   -20160 y_{1} + 1260 y_{2} -  \frac{2240}{27}y_{3}  +  \frac{161000}{9}hy_{0}^{(1)}
     \right. \nonumber \\
    & & \ \ \ \ \ \ \ \ \ \left. +\frac{23800}{3}h^{2}y_{0}^{(2)}+\frac{6160}{3}h^{3}y_{0}^{(3)}+280h^{4}y_{0}^{(4)}
     \right ], \\ \label{e:3.2}
     U_{1} -\frac{30666}{8867}U_{2} +  U_{5} & = & \frac{1}{h^{7}} \ \left[  \ \frac{-957600}{8867} y_{1}
   +  \frac{1048320}{8867} y_{2} -  \frac{90720}{8867} y_{3}- \frac{866880}{8867}hy_{0}^{(1)}
          \right. \nonumber \\
    & & \ \ \ \ \ \ \ \ \ \left.  -  \frac{1209600}{8867} \ h^{2}y_{0}^{(2)}
    -  \frac{829920}{8867} \ h^{3}y_{0}^{(3)} -  \frac{352800}{8867} \ h^{4}y_{0}^{(4)} \right. \nonumber \\
    & & \ \ \ \ \ \ \ \ \ \left. -  \frac{87864}{8867} \ h^{5}y_{0}^{(5)}
     \right ], \nonumber \\
     & &  \\ \label{e:3.3}
     U_{2} -\frac{278026}{94221}U_{3} +  U_{6} & = & \frac{1}{h^{7}} \ \left[  \ \frac{-67340}{10469} y_{2}
   +  \frac{80640}{10469} y_{3} -  \frac{700}{551} y_{4} - \frac{54040}{10469}hy_{0}^{(1)}
         \right. \nonumber \\
       & & \ \ \ \ \ \ \ \ \ \left. -  \frac{4200}{361} \ h^{2}y_{0}^{(2)}
    -  \frac{20720}{1653}h^{3}y_{0}^{(3)}-  \frac{85400}{10469}h^{4}y_{0}^{(4)} \right. \nonumber \\
    & & \ \ \ \ \ \ \ \ \ \left. -  \frac{95536}{31407}h^{5}y_{0}^{(5)}
     \right ], \\ \label{e:3.4}
     U_{3}  +   U_{7}  & = & \frac{1}{h^{7}} \ \left[ \ \frac{10808537040}{487056529}y_{3}
   -  \frac{13373418240}{487056529}y_{4} +  \frac{2564881200}{487056529} y_{5}
    \right. \nonumber \\
   & & \left.
    +  \frac{8243655840}{487056529}hy_{0}^{(1)} +   \frac{26287914240}{487056529} \ h^{2}y_{0}^{(2)} \right. \nonumber \\
   & & \left.
    + \frac{40576352880}{487056529}h^{3}y_{0}^{(3)} +   \frac{39377200800}{487056529} \ h^{4}y_{0}^{(4)} \right. \nonumber \\
   & & \left.
    + \frac{25438766892}{487056529}h^{5}y_{0}^{(5)}+ \frac{9474762304}{487056529}h^{6}y_{0}^{(6)}
    \right], \\ \label{e:3.5}
     U_{4}  +   U_{8}  & = & \frac{1}{h^{7}} \ \left[ \ \frac{2645350155}{436783036}y_{4}
   -  \frac{869117760}{109195759}y_{5} +  \frac{831120885}{436783036} y_{6}
    \right. \nonumber \\
   & & \left.
    +  \frac{31279815}{7530742}hy_{0}^{(1)} +   \frac{3666455415}{218391518} \ h^{2}y_{0}^{(2)} \right. \nonumber \\
   & & \left.
    + \frac{3572264955}{109195759}h^{3}y_{0}^{(3)} +   \frac{8717751945}{218391518} \ h^{4}y_{0}^{(4)} \right. \nonumber \\
   & & \left.
    + \frac{3525702999}{109195759}h^{5}y_{0}^{(5)}+ \frac{1634628387}{109195759}h^{6}y_{0}^{(6)}
    \right]
   \end{eqnarray}
and
\begin{eqnarray}\label{e:3.6}
  U_{5}  +   U_{9}  & = & \frac{1}{h^{7}} \ \left[ \ \frac{132838307280}{61865369749}y_{5}
   -  \frac{183300929280}{61865369749}y_{6} +  \frac{50462622000}{61865369749} y_{7}
    \right. \nonumber \\
   & & \left.
    +  \frac{82375685280}{61865369749}hy_{0}^{(1)} +   \frac{402603647040}{61865369749} \ h^{2}y_{0}^{(2)} \right. \nonumber \\
   & & \left.
    + \frac{946588828080}{61865369749}h^{3}y_{0}^{(3)} +   \frac{1390554453120}{61865369749} \ h^{4}y_{0}^{(4)} \right. \nonumber \\
   & & \left.
    + \frac{1350858565644}{61865369749}h^{5}y_{0}^{(5)}+ \frac{749461929944}{61865369749}h^{6}y_{0}^{(6)}
    \right]\ .
    \end{eqnarray}
    Basically, one does a power expansion in $h$, and postulates that the low orders in $h$ vanish.
    Our calculations were done using Mathematica.

    The local truncation errors associated with the linear equations
    $(\ref{e:3.1})-(\ref{e:3.6})$ and $(\ref{e:2.8})$ are calculated as
\begin{equation}\label{e:3.12}
\tilde{t}_{i} = \left\{ \begin{array}{ll}
   - 5.778 h^{9}y^{(9)}(t_{1})+O(h^{10}), &  i=1,  \\
    - 6.472 h^{9}y^{(9)}(t_{2})+O(h^{10}), &  i=2,  \\
    - 7.230 h^{9}y^{(9)}(t_{3})+O(h^{10}), &  i=3,  \\
    - 19.288 h^{9}y^{(9)}(t_{4})+O(h^{10}), &  i=4,  \\
    - 25.620 h^{9}y^{(9)}(t_{5})+O(h^{10}), &  i=5,  \\
    - 33.020 h^{9}y^{(9)}(t_{6})+O(h^{10}), &  i=6,  \\
   \frac{1}{2}(-100+25 \alpha +13 \beta+5 \gamma+\delta) h^{9}y^{(9)}(t_{i})+O(h^{10}), &  i=7,8,\ldots,n  \\
        \end{array} \right.
\end{equation}
and
\begin{equation}\label{e:3.13}
\|\tilde{T}\|=ch^{9}R=O(h^{9}), \ \ \ \ \ \ \ \ R= \max_{t \in [a, \
b]} |y^{(9)}(t)|,
\end{equation}\
 where $c$ is a constant which
depends only on the values of $\alpha, \ \beta, \ \gamma$ and
$\delta$ and is independent of $h$. Moreover, $\alpha$, $\beta$,
$\gamma$
and $\delta$ are taken such that $\alpha+\beta+\gamma+\delta=60$. \\
In general, the solution of the system of linear equations
$(\ref{e:3.1})-(\ref{e:3.6})$ and $(\ref{e:2.8})$ is second order
convergent. \\
The local truncation error of the system
(\ref{e:2.8}) can be expressed in the following form
\begin{equation}\label{e:4321}
\tilde{t}_{i} = \left\{ \begin{array}{ll}
      2(-60+ \alpha + \beta+ \gamma+\delta) h^{7}y^{(7)}(t_{i})+(-60+ \alpha + \beta+ \gamma+\delta) h^{8}y^{(8)}(t_{i})
     & \\ +\frac{1}{2}(-100+25 \alpha +13 \beta+5\gamma+\delta) h^{9}y^{(9)}(t_{i}) & \\
      +\frac{1}{6}(-120+37 \alpha +19 \beta+7 \gamma+\delta) h^{10}y^{(10)}(t_{i})& \\
      + \frac{1}{24}(-228+337 \alpha +97 \beta+17\gamma+\delta) h^{11}y^{(11)}(t_{i}) & \\
      +\frac{1}{120}(-380+781 \alpha +211 \beta+31\gamma+\delta) h^{12}y^{(12)}(t_{i})& \\
      +O(h^{13}), &  \\
      i=7,8,\ldots,n,  & \\
 \end{array} \right.
\end{equation}
therefore, the order of the truncation error $\tilde{t}_{i}$ can be
improved to be of order $h^{12}$, if
$\alpha=\frac{151}{15}-\frac{\delta}{5}, \
\beta=-\frac{301}{6}+\delta, \
\gamma=\frac{1001}{10}-\frac{9\delta}{5}.$ Correspondingly, the end
conditions with local truncation error of $O(h^{12})$ can be
determined as
\begin{eqnarray}\label{e:3.14}
  & & U_{0} -\frac{24407}{109} U_{1}-\frac{59362}{109} U_{2} -\frac{10662}{109}U_{3}-\frac{907}{109}U_{4}+
   U_{5} \nonumber \\
   & = & \frac{1}{h^{7}} \ \left[  \ \frac{66633336}{545} y_{0}
   - \frac{19958400}{109} y_{1} + \frac{9979200}{109} y_{2} -  \frac{4435200}{109}y_{3} + \frac{1247400}{109} y_{4}
    - \frac{798336}{545}y_{5}
     \right. \nonumber \\
    & & \ \ \ \ \ \ \ \ \ \left. +  \frac{9114336}{109}hy_{0}^{(1)}+ \frac{1995840}{109}h^{2}y_{0}^{(2)}
    - \frac{80}{109}h^{7}y_{0}^{(7)}
     \right ], \\ \label{e:3.15}
    & & U_{1}+\frac{202055040421}{554613069}U_{2}+\frac{77878525838}{184871023}U_{3} +\frac{18661788874}{184871023}U_{4}
     -\frac{2434662535}{554613069} U_{5}+  U_{6} \nonumber \\
      & = & \frac{1}{h^{7}} \ \left[  \ \frac{-4011644165760}{184871023} y_{1}
   +  \frac{8861887188000}{184871023} y_{2} -  \frac{73815832428800}{1663839207} y_{3}
             \right. \nonumber \\
    & & \ \ \ \ \ \ \ \ \ \left. +   \frac{4618760731200}{184871023}y_{4} - \frac{1470208924800}{184871023}y_{5}
    + \frac{1826678971040}{1663839207} y_{6}  \right. \nonumber \\
    & & \ \ \ \ \ \ \ \ \ \left. -
     \frac{4345911046400}{554613069}hy_{0}^{(1)} -  \frac{1035868310400}{184871023} \ h^{2}y_{0}^{(2)}
    -  \frac{183470425600}{184871023}h^{3}y_{0}^{(3)}
     \right ], \nonumber \\
     & &  \\ \label{e:3.16}
   & &  \left[U_{2}-\frac{13173366154319505819}{604803696004634}U_{3}-\frac{5923535526089565973}{302401848002317}U_{4}
      \right. \nonumber \\
     & & \left. - \frac{2639790737228529743}{604803696004634} U_{5} +
     U_{7} \right] \nonumber \\
     & = & \frac{1}{h^{7}} \ \left[  \ \frac{25352931909798309915}{43200264000331} y_{2}
   -  \frac{198185856313975120000}{129600792000993} y_{3} \right. \nonumber \\
    & & \ \ \ \ \ \ \ \ \ \left. +  \frac{72535593878062755750}{43200264000331} y_{4}
    -   \frac{44672825515677652800}{43200264000331}y_{5}
\right. \nonumber \\
    & & \ \ \ \ \ \ \ \ \ \left.
    +\frac{44958164899589796925}{129600792000993}y_{6}
    -\frac{2139803134054971840}{43200264000331} y_{7}  \right. \nonumber \\
    & & \ \ \ \ \ \ \ \ \ \left. +
     \frac{5764036699720950200}{43200264000331}hy_{0}^{(1)}
     +  \frac{7374675959642702700}{43200264000331} \ h^{2}y_{0}^{(2)}
    \right. \nonumber \\
    & & \ \ \ \ \ \ \ \ \ \left. +  \frac{3273172503578299200}{43200264000331}h^{3}y_{0}^{(3)}
    +  \frac{521467194925746900}{43200264000331}h^{4}y_{0}^{(4)}
     \right ], \\ \label{e:3.17}
    & & U_{3} - \frac{3169885805313999875741}{17618985121607404312}U_{4}
    - \frac{544760103861334609083}{17618985121607404312}U_{6} +   U_{7} \nonumber \\
      & = & \frac{1}{h^{7}} \ \left[ \ \frac{5366584014500607349360}{19821358261808329851}y_{3}
   -  \frac{1627194491533397967735555}{1127615047782873875968}y_{4} \right. \nonumber \\
    & & \ \ \ \ \ \ \ \ \ \left.
   +  \frac{5581229831865256388400}{2202373140200925539} y_{5}
    -   \frac{324347788172453021646845}{158570866094466638808}y_{6}
    \right. \nonumber \\
   & & \ \ \ \ \ \ \ \ \  \left.
    + \frac{1793040053953186573920}{2202373140200925539}y_{7} - \frac{147257531382013448691645}{1127615047782873875968}y_{8}
    \right. \nonumber \\
   & & \ \ \ \ \ \ \ \ \  \left. -  \frac{78317811652764229087465}{845711285837155406976}hy_{0}^{(1)} -   \frac{42412130734120984954425}{140951880972859234496} \
   h^{2}y_{0}^{(2)} \right. \nonumber \\
   & & \ \ \ \ \ \ \ \ \ \left.   -\frac{13509365844145615609375}{35237970243214808624}h^{3}y_{0}^{(3)}  -\frac{7826257773953806554675}{35237970243214808624}h^{4}y_{0}^{(4)} \right. \nonumber \\
   & & \ \ \ \ \ \ \ \ \ \left. -\frac{447768854035545682017}{8809492560803702156}h^{5}y_{0}^{(5)}
    \right], \\ \label{e:3.18}
    & & U_{4} - \frac{6169811365491003355386625}{364845537886699795641421}U_{7}
    +   U_{9} \nonumber \\
      & = & \frac{1}{h^{7}} \ \left[ \ \frac{413182203198678792199193360481}{373601830795980590736815104}y_{4}
   -
   \frac{1011762526223941981900336800}{364845537886699795641421}y_{5} \right. \nonumber \\
   & & \ \ \ \ \ \ \ \ \  \left.
   +  \frac{998387082478934194463004566965}{354629862825872201363461212} y_{6}
    -   \frac{566429213407879867786917120}{364845537886699795641421}y_{7}
    \right. \nonumber \\
   & & \ \ \ \ \ \ \ \ \ \left.
    + \frac{173230355267937275186019127455}{373601830795980590736815104}y_{8}
    -
    \frac{5254626822196644195075230752}{88657465706468050340865303}y_{9} \right. \nonumber \\
   & & \ \ \ \ \ \ \ \ \  \left.
    +  \frac{1828802733123508354716945025585}{7565437073618606962420505856}hy_{0}^{(1)} \right. \nonumber \\
   & & \ \ \ \ \ \ \ \ \ \left. +   \frac{54937023892836663800655898465}{74170951702143205513926528} \
   h^{2}y_{0}^{(2)} \right. \nonumber \\
    & &  \ \ \ \ \ \ \ \ \ \left. + \frac{319183920456421230207708911935}{315226544734108623434187744}h^{3}y_{0}^{(3)} \right. \nonumber \\
   & & \ \ \ \ \ \ \ \ \ \left. + \frac{81879551659139637198985554365}{105075514911369541144729248}h^{4}y_{0}^{(4)} \right. \nonumber \\
   & & \ \ \ \ \ \ \ \ \ \left. +
   \frac{2991554077139003376141526763}{8756292909280795095394104}h^{5}y_{0}^{(5)} \right. \nonumber \\
   & & \ \ \ \ \ \ \ \ \ \left. + \frac{303485670688565607390252013}{4378146454640397547697052}h^{6}y_{0}^{(6)}
    \right]
   \end{eqnarray}
and
\begin{eqnarray}\label{e:3.180}
 & & U_{5}  +   U_{10} \nonumber \\
    & = & \frac{1}{h^{7}} \ \left[ \ \frac{19038680213948167651954555270266}{43087137994818537402205515625}y_{5}
   - \frac{3612553213748861716357961962885}{2680364680381671574716400716}
   y_{6} \right. \nonumber \\
   & & \ \ \ \ \ \ \ \ \  \left.
   + \frac{4872382360659412438888663650}{2757576831668386393741153}  y_{7}
   -   \frac{3637650046079073899112436800}{2757576831668386393741153}
   y_{8} \right. \nonumber \\
   & & \ \ \ \ \ \ \ \ \  \left.
    + \frac{396380487483749310137552942650}{670091170095417893679100179} y_{9} \right. \nonumber \\
   & & \ \ \ \ \ \ \ \ \left.
     - \frac{25866940054548307411938055846689}{172348551979274149608822062500} y_{10} \right. \nonumber \\
   & & \ \ \ \ \ \ \ \ \left.  +
   \frac{38692597856846972600037888063421}{698011635516060305915729353125}hy_{0}^{(1)} \right. \nonumber \\
   & & \ \ \ \ \ \ \ \ \left.
   +
   \frac{9138473244848295227157642242648}{46534109034404020394381956875}h^{2}y_{0}^{(2)} \right. \nonumber \\
   & & \ \ \ \ \ \ \ \ \left.
   +   \frac{2879180812841847947352594547754}{9306821806880804078876391375} \
   h^{3}y_{0}^{(3)} \right. \nonumber \\
   & & \ \ \ \ \ \ \ \ \left. +   \frac{1170428496335992779989052992}{4278998531899220266150065} \
   h^{4}y_{0}^{(4)} \right. \nonumber \\
   & & \ \ \ \ \ \ \ \ \left. +   \frac{2186365487813281315497916274}{15909097105779152271583575} \
   h^{5}y_{0}^{(5)} \right. \nonumber \\
   & & \ \ \ \ \ \ \ \ \left. +   \frac{1319003601532979667100927096}{41363652475025795906117295} \
   h^{6}y_{0}^{(6)}    \right] \ .
    \end{eqnarray}
It turns out that if $\alpha$, $\beta$, $\gamma$ and $\delta$ are
chosen as $\alpha=\frac{151}{15}-\frac{\delta}{5}, \
\beta=-\frac{301}{6}+\delta, \
\gamma=\frac{1001}{10}-\frac{9\delta}{5},$ then the order of
truncation error of eq. $(\ref{e:2.8})$ is $O(h^{12})$ and the order
of convergence can be improved up to
five based on the improved order of the end conditions.\\ \\
To illustrate the implementation of the method, three examples are
discussed in the following section. We choose a driving force at the top of the cascade
for which the IVP can be solved analytically, and then
investigate the error terms of our numerical method by comparing with the exact solution.
\newpage
 \newsection{Numerical Examples}
   \ \ \\
{\large {\bf Example 1 }} \\
  Consider the following initial value problem
 \begin{equation}
\left.
\begin{array}{rl}\label{e:4.1}
    y^{(7)}(t) \ + \  y(t) \ & =  \ - (t^{2}-43) \ cos(t)+(-1+t^{2}-14t) sin(t), \ \ \ \ \ t \in [-1, \ 1],  \\
 y(-1) \ & =  0, \\
   y^{(1)}(-1) \ & =  \ 2sin(1), \\
   y^{(2)}(-1) \ & =  -4cos(1)-2sin(1)\\
   y^{(3)}(-1) \ & = 6cos(1)-6sin(1)\\
   y^{(4)}(-1) \ & = 8cos(1)+12sin(1)\\
   y^{(5)}(-1) \ & = -20cos(1)+10sin(1)\\
   y^{(6)}(-1) \ & = -12cos(1)-30sin(1).
\end{array}
\right\}
\end{equation}
This basically corresponds to a periodic forcing at the top of the cascade
whose amplitude is a particular quadratic function of $t$.

The analytic solution of the above problem is
$$ y(t) \ = \ (t^{2}-1) \  sin(t) \ ,$$
meaning the velocity of the driven particle oscillates with increasing
amplitude. Since we can solve this example analytically, we can easily
determine the numerical error of our method.
The observed maximum errors (in absolute values) associated with $
y_{i}$, for the problem $(\ref{e:4.1})$, corresponding to the
different values of $\alpha, \ \beta, \ \gamma$ and $\delta$, are
summarized in Table 1.
\par \noindent
\begin{table}[htp]
\caption{ Maximum absolute errors for problem $(\ref{e:4.1})$ in
 $ y_{i}.$ }
 \begin{center}
\begin{tabular}{|c|c|c|c|} \hline
               &                              &
               &                                           \\
$n$              &  $\alpha=\frac{1}{2}, \ \beta=\frac{19}{2} $
                &
                  $\alpha=0, \ \beta=0 $  &
                  $\alpha=10, \ \beta=10 $
                    \\
               &               &
               &
                                                       \\
               &  $\gamma=\frac{49}{2}, \ \delta=\frac{51}{2} $             &
              $\gamma=0, \ \delta=60 $ &
              $\gamma=10, \ \delta=30 $
                                                      \\
               &                &
               &
                                                            \\ \hline
               &                              &
               &                                            \\
$12 $        & $ 2.88 \times {10}^{-1} $
               & $ 3.04 \times {10}^{-1} $
               & $ 2.76 \times {10}^{-1} $
                                                                      \\ \hline
               &                              &
               &                                                                                 \\
$ 24$   & $ 3.09 \times {10}^{-2} $
               & $ 3.56 \times {10}^{-2} $
               & $ 2.73 \times {10}^{-2} $
                                                                                     \\ \hline
               &                              &
               &                                                                                 \\
$48 $  & $ 2.5 \times {10}^{-3} $
               & $ 3.9 \times {10}^{-3} $
               & $ 1.4 \times {10}^{-3} $
                                                                  \\ \hline
               &                              &
               &                                                                                \\
$96$ & $ 1.70 \times {10}^{-4} $
               & $ 7.37 \times {10}^{-4} $
               & $ 3.19 \times {10}^{-4} $
                                               \\ \hline
 \end{tabular}
\end{center}
\end{table}
The observed maximum errors (in absolute values) associated with $
y_{i}$, for the problem $(\ref{e:4.1})$, corresponding to the
improved end conditions are summarized in Table 2. A significant improvement
is obtained.
\par \noindent
\begin{table}[htp]
\caption{ Maximum absolute errors for problem $(\ref{e:4.1})$ in
 $ y_{i}.$ }
 \begin{center}
\begin{tabular}{|c|c|} \hline
               &                                          \\
$n$        & $|y(t_{i})-y_{i} |$
                                                            \\ \hline
               &                                            \\
$10$       & $ 2.25 \times {10}^{-1} $
                                                                                     \\ \hline
               &                                                                               \\
$ 20$   & $ 2.08 \times {10}^{-6} $
                                                                                                    \\ \hline
               &                                                                                \\
$40 $  & $ 7.50 \times {10}^{-7}$
                                                                  \\ \hline
 \end{tabular}
\end{center}
\end{table}
\par \noindent
{\large {\bf Example 2 }} \\
 \ \ \\
 Consider now an exponential forcing at the top of
 the cascade, namely the following IVP:
 \begin{equation}
\left.
\begin{array}{ll}\label{e:4.2}
    y^{(7)}(t) \ -  y(t) \ =  \ -7 e^{t}(5+2t), \ \ \ \ \ 0\leq t \leq 1 \\
 y(0) \  \ \ \  = 0,  \ \ \ \ \ \ \ \ \ \ \ \ \ \ \ \  y^{(1)}(0) \ = 1, \\
   y^{(2)}(0) \ =  0, \ \ \ \ \ \ \ \ \ \ \ \ \ y^{(3)}(0) \ =  -3 \\
   y^{(4)}(0) \ =  -8, \ \ \ \ \ \ \ \ \ \ \ \ \ y^{(5)}(0) \ = -15, \\
   y^{(6)}(0) \ =  -24.
\end{array}
\right\}
\end{equation}
The corresponding analytic solution is
$$ y(t) \ = \ t(1-t) \  e^{t} \ .$$
In this case, for a limited amount of time, velocity differences grow exponentially,
consistent with the chaotic nature of the velocity field in turbulent flows.
In this case the observed maximum errors of our method are shown in table 3.

\par \noindent
\begin{table}[htp]
\caption{ Maximum absolute errors for problem $(\ref{e:4.2})$ in
 $ y_{i}.$ }
 \begin{center}
\begin{tabular}{|c|c|c|c|} \hline
               &                              &
               &                                           \\
$n$              &  $\alpha=\frac{1}{2}, \ \beta=\frac{19}{2} $
                &
                  $\alpha=0, \ \beta=0 $  &
                  $\alpha=10, \ \beta=10 $
                    \\
               &               &
               &
                                                       \\
               &  $\gamma=\frac{49}{2}, \ \delta=\frac{51}{2} $             &
              $\gamma=0, \ \delta=60 $ &
              $\gamma=10, \ \delta=30 $
                                                      \\
               &                &
               &
                                                            \\ \hline
               &                              &
               &                                            \\
$10 $        & $ 1.5 \times {10}^{-3} $
               & $ 1.6 \times {10}^{-3} $
               & $ 1.5 \times {10}^{-3} $
                                                                      \\ \hline
               &                              &
               &                                                                                 \\
$ 20$   & $ 1.75 \times {10}^{-4} $
               & $ 1.94 \times {10}^{-4} $
               & $ 1.60 \times {10}^{-4} $
                                                                                     \\ \hline
               &                              &
               &                                                                                 \\
$40$  & $ 1.81 \times {10}^{-5} $
               & $ 2.62 \times {10}^{-5} $
               & $ 1.32 \times {10}^{-5} $
                                                                  \\ \hline
 \end{tabular}
\end{center}
\end{table}
The observed maximum errors (in absolute values) associated with $
y_{i}$, for the problem $(\ref{e:4.2})$, corresponding to the
improved end conditions are summarized in Table 4.
\par \noindent
\begin{table}[htp]
\caption{ Maximum absolute errors for problem $(\ref{e:4.2})$ in
 $ y_{i}.$ }
 \begin{center}
\begin{tabular}{|c|c|} \hline
               &                                           \\
$n$            & $|y(t_{i})-y_{i} |$
                                                            \\ \hline
               &                                            \\
$10 $        & $ 1.82 \times {10}^{-1} $
                                                                                     \\ \hline
               &                                                                               \\
$ 12$  & $ 2.15 \times {10}^{-8} $
                                                                                                    \\ \hline
               &                                                                                \\
$15$   & $ 3.65 \times {10}^{-9}$
                                                                  \\ \hline
 \end{tabular}
\end{center}
\end{table}
\par \noindent
{\large {\bf Example 3 }} \\
 \ \ \\
 Consider the following problem:
 \begin{equation}\label{e:4.3}
\left.
\begin{array}{ll}
    y^{(7)}(t) \  \ =  \ 7(-6e^{t}+e^{t}(1-t))+(-7e^{t}+e^{t}(1-t))t, \ \ \ \ \ 0\leq t \leq 1 \\
 y(0) \  \ \ \  = 0,  \ \ \ \ \ \ \ \ \ \ \ \ \ \ \ \  y^{(1)}(0) \ = 1, \\
   y^{(2)}(0) \ =  0, \ \ \ \ \ \ \ \ \ \ \ \ \ y^{(3)}(0) \ =  -3 \\
   y^{(4)}(0) \ =  -8, \ \ \ \ \ \ \ \ \ \ \ \ \ y^{(5)}(0) \ = -15, \\
   y^{(6)}(0) \ =  -24.
\end{array}
\right\}
\end{equation}
The corresponding analytic solution is again
$$ y(t) \ = \ t(1-t) \  e^{t} \ $$
and we get in this case the results shown in table 5 and 6.
\par \noindent
\begin{table}[htp]
\caption{ Maximum absolute errors for problem $(\ref{e:4.3})$ in
 $ y_{i}.$ }
 \begin{center}
\begin{tabular}{|c|c|c|c|} \hline
               &                              &
               &                                           \\
$n$              &  $\alpha=\frac{1}{2}, \ \beta=\frac{19}{2} $
                &
                  $\alpha=0, \ \beta=0 $  &
                  $\alpha=10, \ \beta=10 $
                    \\
               &               &
               &
                                                       \\
               &  $\gamma=\frac{49}{2}, \ \delta=\frac{51}{2} $             &
              $\gamma=0, \ \delta=60 $ &
              $\gamma=10, \ \delta=30 $
                                                      \\
               &                &
               &
                                                            \\ \hline
               &                              &
               &                                            \\
$9$        & $ 2.0 \times {10}^{-3} $
               & $ 2.22 \times {10}^{-3} $
               & $ 1.5 \times {10}^{-3} $
                                                                      \\ \hline
               &                              &
               &                                                                                 \\
$ 18$   & $ 2.26 \times {10}^{-4} $
               & $ 2.66 \times {10}^{-4} $
               & $ 1.60 \times {10}^{-4} $
                                                                                     \\ \hline
               &                              &
               &                                                                                 \\
$36$  & $ 2.16 \times {10}^{-5} $
               & $ 3.46 \times {10}^{-5} $
               & $ 1.32 \times {10}^{-5} $
                                                                  \\ \hline
 \end{tabular}
\end{center}
\end{table}
\par \noindent
\begin{table}[htp]
\caption{ Maximum absolute errors for problem $(\ref{e:4.3})$ in
 $ y_{i}.$ }
 \begin{center}
\begin{tabular}{|c|c|} \hline
               &                                           \\
$n$            & $|y(t_{i})-y_{i} |$
                                                            \\ \hline
               &                                            \\
$10 $        & $ 1.82 \times {10}^{-1} $
                                                                                     \\ \hline
               &                                                                               \\
$ 12$  & $ 2.33 \times {10}^{-8} $
                                                                                                    \\ \hline
               &                                                                                \\
$15$   & $ 1.67 \times {10}^{-8}$
                                                                  \\ \hline
 \end{tabular}
\end{center}
\end{table}

\par \noindent
{\large {\bf Conclusion}} \\

In this paper we showed that hierarchical cascade models, motivated
by turbulent flows, can lead to initial value problems of $N$-th
order, where $N$ is the number of cascade steps. As an example we
considered $N=7$ and designed the optimum strategy to numerically
treat the corresponding initial value problem. A non-polynomial
spline method was developed for this. The numerical algorithm
depends on some parameters $\alpha, \beta, \gamma, \delta$ for which
we derived explicit formulas. The method is observed to be 
second-order convergent for arbitrary choices of the parameters $\alpha$,
$\beta$, $\gamma$ and $\delta$ such that
$\alpha+\beta+\gamma+\delta=60$ but if $\alpha$, $\beta$, $\gamma$
and $\delta$ are chosen as $\alpha=\frac{151}{15}-\frac{\delta}{5},
\ \beta=-\frac{301}{6}+\delta, \
\gamma=\frac{1001}{10}-\frac{9\delta}{5},$ then the method is of
order five, due to the use of improved order end conditions. This is
the optimum choice to numerically deal with cascade-like models of
this type. In turbulent flows, driving forces and velocity
differences often behave chaotically, corresponding to exponential
growth for a limited amount of time. Consequently, we tested our
method for these and other types of driving forces. We found that
the maximum numerical errors observed are indeed very small if
improved order end conditions are used.

\end{document}